%% file: OCIN-ArXives.tex
\documentclass[11pt]{amsart}

\usepackage[usenames]{color}

\usepackage[colorlinks=true,
linkcolor=webgreen,
filecolor=webbrown,
citecolor=webgreen]{hyperref}

\definecolor{webgreen}{rgb}{0,.5,0}
\definecolor{webbrown}{rgb}{.6,0,0}

\newcommand{\seqnum}[1]{\href{http://oeis.org/#1}{\underline{#1}}}

\DeclareMathOperator{\dom}{Dom}
\DeclareMathOperator{\im}{Im}

\include{cdmstdcmds}

\begin{document}

 \title[Combinatorial results for contraction mappings]{Combinatorial results for order-preserving partial injective contraction mappings}

\author[B.M.~Ahmed]{Bayo Musa Ahmed}
\address{Department of Mathematics,
Faculty of Physical Sciences,
Univesity of Ilorin, NIGERIA}
\email{ahmed.bm@unilorin.edu.ng}

\author[N.~AlDhamri]{Nadia AlDhamri}
\address{Department of Mathematics and Statistics,
Sultan Quaboos University,
Al-Khod, PC 123, OMAN }
\email{nadias@squ.edu.om}

\author[F.~Al-Kharousi]{Fatma Al-Kharousi}
\address{Department of Mathematics and Statistics,
Sultan Quaboos University,
Al-Khod, PC 123, OMAN }
\email{fatma9@squ.edu.om}

\author[G.~Klein]{Georg Klein}
\address{Departamento de Matem\'atica,
Universidade Federal da Bahia,
Av. Adhemar de Barros, S/N, Ondina, CEP: 40.170.110, Salvador, BA, BRAZIL}
\email{georgklein53@gmail.com}

\author[A.~Umar]{Abdullahi Umar}
\address{Department of Mathematical Sciences,
Khalifa University of Science and Technology, 
SAN Campus, PO. Box 127788
Abu Dhabi, UAE}
\email{abdullahi.umar@ku.ac.ae}


\thanks{Financial support from The Reasearch Council of Oman grant RC/SCI/DOMS/13/01 is gratefully acknowledged}

\begin{abstract} 
Let $ \mathcal{I}_n$ be the symmetric inverse semigroup on \\ $X_n = \{1,
2, \ldots , n\}$. Let $\mathcal{OCI}_n$ be the subsemigroup of  $\mathcal{I}_n$ consisting of 
all order-preserving injective partial contraction mappings,
and let $\mathcal{ODCI}_n$ be the subsemigroup of  $\mathcal{I}_n$ consisting of
all order-preserving and order-decreasing injective partial contraction mappings of $X_n$.
In this paper, 
we investigate the cardinalities of some equivalences on $\mathcal{OCI}_n$ and $\mathcal{ODCI}_n$ which lead naturally
to obtaining the order of these semigroups. Then, we relate the formulae obtained to Fibonacci numbers. 
Similar results about $\mathcal{ORCI}_n$, the semigroup of 
order-preserving or order-reversing injective partial contraction mappings, are deduced. 
\end{abstract}

\maketitle


\section{Introduction and Preliminaries}\label{introduction} 

Let $X_n=\{1,2, \ldots, n\}$ and $\mathcal{I}_n$ be the partial
one-to-one transformation semigroup on $X_n$ under composition of
mappings. Then $\mathcal{I}_n$ is an {\em inverse} semigroup, that is,
for all $\alpha \in \mathcal{I}_n$ there exists a unique $\alpha' \in
\mathcal{I}_n$ such that $\alpha = \alpha\alpha'\alpha$ and $\alpha' =
\alpha'\alpha\alpha'$. The importance of $\mathcal{I}_n$, 
commonly known as the {\em symmetric inverse semigroup or monoid}, to
inverse semigroup theory is similar to the importance of the symmetric
group $\mathcal{S}_n$ to group theory. Every finite inverse semigroup
$S$ is embeddable in $\mathcal{I}_n$, an analogue to Cayley's theorem
for finite groups. Thus, just as the study of symmetric, alternating, and
dihedral groups has significantly contributed to group theory,
the study of various subsemigroups of $\mathcal{I}_n$ lead to significant 
contributions to the theory of semigroups. For instance,  see 
Borwein et al.~\cite{Bor},  Fernandes~\cite{Fer1}, Fernandes et al.~\cite{Fer2},
Garba~\cite{Gar},  Laradji and Umar~\cite{Lar},
and Umar~\cite{Uma1, Uma2}.

We shall denote the domain of $\alpha \in \mathcal{I}_n$ by $\dom \alpha$. A transformation $\alpha \in \mathcal{I}_n$ is said to be
{\em order-preserving} if for all $x,y$ in 
$\dom \alpha$, $ \displaystyle{ x \leq y }$ implies $x\displaystyle{ \alpha \leq y\alpha }$.
A transformation $\alpha \in \mathcal{ I}_n$ is said to be
{\em order-reversing} if for all $x,y$ in 
$\dom \alpha$, $ \displaystyle{ x \leq y }$ implies $x\displaystyle{ \alpha \geq y\alpha }$.
A transformation $\alpha \in \mathcal{ I}_n$ is an {\em isometry} if it is  {\em distance-preserving}, i.e, 
for all $x,y$ in 
$\dom \alpha$, $ | x-y| = |x\alpha -y\alpha |$. 
A transformation $\alpha \in \mathcal{ I}_n$ is a  {\em  contraction} if for all $x,y$ in 
$\dom \alpha$, $ | x\alpha-y\alpha  | \leq | x -y|$. A transformation is
said to be {\em order-decreasing} if for all $x$ in 
$\dom \alpha$, $\displaystyle{ x \alpha \leq x}$. 

Analogous to Al-Kharousi et al.~\cite{Kha2,Kha3}, we investigate the
combinatorial properties of $\mathcal{ OCI}_n$ and $\mathcal{ ODCI}_n$, thereby
complementing the results in Al-Kharousi et al.~\cite{Kha1} which deals primarily
with the algebraic and rank properties of $\mathcal{OCI}_n$.
In the present section, we introduce basic definitions and terminology, and we
quote some elementary results from Al-Kharousi et al.~\cite{Kha1,Kha2} that will be needed 
subsequently. 
In Section~\ref{CombOCIn},
we thoroughly  investigate the combinatorial properties of $\mathcal{ OCI}_n$, the semigroup of 
order-preserving partial injective contraction mappings of the finite chain $X_n$. 
We obtain a formula for $F(n;p)$, the number of transformations with height $p$. 
Then we refine this formula obtaining $F(n;p,m)$, the number of transformations with height $p$ and $m$ fixed points.

We recall that the Fibonacci sequence, denoted by $F_n$, is defined recursively such that each number is the sum of the two preceding ones, starting from 0 and 1. 
The order 
of $\mathcal{ OCI}_n$ is obtained and expressed in terms of two consecutive Fibonacci numbers, and it is shown that this is  
Sequence \seqnum{A094864} in The On-Line Encyclopedia of
Integer Sequences~\cite{oeis}. 
In Section~\ref{CombODCIn}, 
we study the combinatorial properties of $\mathcal{ ODCI}_n$, the semigroup of 
order-preserving and order-decreasing partial injective contraction mappings of the finite chain $X_n$. 
The investigations go along the lines of the ones in  Section~\ref{CombOCIn} in the case  of $\mathcal{ OCI}_n$. 
The order of $\mathcal{ ODCI}_n$, as a function of $n$, is shown to be 
Sequence \seqnum{A001519} in The On-Line Encyclopedia of
Integer Sequences~\cite{oeis}.
In Section~\ref{CombORCIn},
we study combinatorial properties of $\mathcal{ ORCI}_n$, the semigroup of 
order-preserving or order-reversing partial injective contraction mappings  of the finite chain $X_n$.

For standard concepts in semigroup and symmetric inverse semigroup
theory, we refer to  Howie~\cite{How} and  Lipscomb~\cite{Lip}.

We define the set of all partial injective contractions of $X_n$ as
\begin{eqnarray*}\mathcal{CI}_n = \{ \alpha \in \mathcal{I}_n : \forall x,y \in \dom\alpha, 
\ \mid x\alpha -y\alpha \mid \leq \mid x- y\mid \} \text{.} \end{eqnarray*}
We define the set of all
order-preserving partial injective contractions of $X_n$ as
\begin{eqnarray*}  \mathcal{OCI}_n = \{ \alpha \in \mathcal{CI}_n : \forall x,y \in 
\dom \alpha, \ x \leq y \implies x\alpha \leq y\alpha\} \text{.} \end{eqnarray*}
We define the set of all order-preserving or order-reversing partial injective contractions of $X_n$ as 
\begin{eqnarray*}  \mathcal{ORCI}_n = \{ \alpha \in \mathcal{CI}_n : \forall x,y \in 
\dom \alpha, \ x \leq y \implies x\alpha \leq y\alpha , \\ \text{ or } \forall x,y \in 
\dom \alpha, \  x \leq y \implies x\alpha \geq y\alpha \}  \text{.} \end{eqnarray*}
We define the set of all order-preserving and order-decreasing injective contractions of $X_n$ as
\begin{eqnarray*}  \mathcal{ODCI}_n = \{ \alpha \in \mathcal{OCI}_n : \forall x \in 
\dom \alpha, \ x\alpha \leq x\} \text{.} \end{eqnarray*} 

From these definitions, we can deduce the following.
\begin{lemma} \label{SubSemi} $\mathcal{ CI}_n$, $\mathcal{ OCI}_n$, $\mathcal{ ORCI}_n$ and $\mathcal{ ODCI}_n$ are
subsemigroups of $\mathcal{ I}_n$.  \qed
\end{lemma} 

Let $\alpha$ be an arbitrary element in $\mathcal{ I}_n$. The {\em height} or
{\em rank} of $\alpha$ is defined as $h(\alpha)= | \im \alpha |$. 
The {\em right 
waist} of $\alpha$ is defined as $w^+(\alpha) = max(\im \alpha)$. 
Similarly, the {\em left
waist} of $\alpha$ is defined as
$w^-(\alpha) = min(\im \alpha)$. 
The {\em right 
shoulder} of $\alpha$ is defined as $\varpi^+(\alpha) = max(\dom \alpha)$.
Similarly, the {\em left
shoulder} of $\alpha$ is defined as 
$\varpi^-(\alpha) = min(\dom \alpha)$.
The {\em fix} of $\alpha$
is denoted by $f(\alpha)$, and defined as $f(\alpha)=|F(\alpha)|$, where
$F(\alpha) = \{x \in X_n: x\alpha = x\}$.

For a given transformation semigroup $S$,
we define $$F(n;p)=|\{ \alpha \in S :  \ h(\alpha )= p   \}| \text{.}$$

The \textit{gap} of an ordered tuple $\mathbf{a}=(a_1,a_2,\ldots,a_p)$ is the ordered $(p-1)$-tuple
$$g(\mathbf{a})=\left(a_{2}-a_{1}, a_{3}-a_{2}, \dots, a_{p}-a_{p-1}\right) \text{.}$$
Accordingly, for $\displaystyle{ \alpha= \left( \begin{array}{cccc}
           a_{1} & a_{2}& \cdots & a_{p}  \\
            a_{1}\alpha & a_{2}\alpha& \cdots & a_{p}\alpha \\
          \end{array} \right)}$  \vspace{3 pt}\\ 
with $\displaystyle{1\leq   a_{1}< a_{2}< \dots< a_{p}\leq n}$, 
let the gap of the domain of $\alpha$ 
be
$$ g(\dom \alpha)=\left(a_{2}-a_{1}, a_{3}-a_{2}, \dots, a_{p}-a_{p-1}\right) \text{,}$$ 
and let the gap of the image of $\alpha$ be  
$$g(\im \alpha)=\left(a_{2}\alpha-a_{1}\alpha, a_{3}\alpha-a_{2}\alpha, \dots, a_{p}\alpha-a_{p-1}\alpha\right) \text{.}$$ 

For example, if $\displaystyle{ \alpha= \left( \begin{array}{ccc}
           1&3&5 \\
          3&5&6\\
          \end{array} \right), \beta= \left( \begin{array}{cccc}
          1  &2&3  &5   \\
            3 &4 &5 &6 \\
          \end{array} \right)\in \mathcal{OCI}_6}$, then \vspace{3 pt} \\ $g(\dom \alpha)=\left(2,2\right)$, $g(\im \alpha)=\left(2,1\right)$, 
          $g(\dom \beta)=\left(1,1,2\right)$ and $g(\im \beta)=\left(1,1,1\right)$. \vspace{3 pt} 

For $\displaystyle{ \alpha= \left( \begin{array}{cccc}
           a_{1} & a_{2}& \cdots & a_{p}  \\
            a_{1}\alpha & a_{2}\alpha& \cdots & a_{p}\alpha \\
          \end{array} \right)}$ \vspace{3 pt} 
with $\displaystyle{1\leq   a_{1}< a_{2}< \dots< a_{p}\leq n}$,  let $\displaystyle{ d_{i}= a_{i+1}\alpha-a_{i}\alpha}$
for $\displaystyle{ i=1,2,\dots ,p-1}$. Then the {\em gap of the image set} of $\alpha$ is the ordered $(p-1)$-tuple $(d_1, d_2, \ldots , d_{p-1})$. 
Similarly, we let $\displaystyle{ t_{i}=a_{i+1}-a_{i}}$ for $\displaystyle{ i=1,2,\dots ,p-1}$. Then  
the {\em gap of the domain set} of $\alpha$  is the ordered $(p-1)$-tuple $(t_1, t_2, \ldots , t_{p-1})$. With these  definitions, it is clear that
$$p-1\leq\sum_{i=1}^{p-1}\left| d_{i} \right| \leq n-1 \text{ \ and \  } p-1\leq\sum_{i=1}^{p-1} t_{i}  \leq n-1  \text{.}$$

Here we state some well known identities which we will use in the proofs of the  results in the following sections. We adopt the convention that 
$\binom{n}{k}=0$ if $k>n$ or if $k$ or $n$ are negative. 
\begin{lemma}\label{tupple}  
\begin{itemize}
\item[(i)]  The number of compositions of $n$ into $p$ positive parts is $\displaystyle  \binom{n-1}{p-1}$, \\
(See \cite[p. 151]{Rio}.).
\item[(ii)]
Let $\displaystyle d(n,p)$ be the number of distinct ordered $p$-tuples: $\displaystyle(r_{1},r_{2},\dots,r_{p})$ 
with $r_i\geq 0$ for $1\leq i\leq p$ and $\displaystyle \sum_{i=1}^{p}r_{i}=n$. 
Then, $d(n,p)$ is the number of \vspace{2 pt} compositions of $n$ into $p$ non-negative parts.
That is, $\displaystyle d(n,p)=  \binom{n+p-1}{p-1}$, \vspace{3 pt} \\ 
(See \cite[p. 589]{Epp}.).  
\end{itemize}  \qed
\end{lemma}


The  Vandermonde convolution identity
states that for integers $m,n,r$,
$$ \sum_{k=0}^r\binom{m}{k}\binom{n}{r-k} =\binom{m+n}{r}   \text{,}$$
(See \cite[p. 8]{Rio}.).
Using this identity we can deduce the following result. 
\begin{lemma}\label{VarVander}
For $r,s,t \in \mathbb{N}$ such that $r> s$, 
$$ \sum_{i=0}^{r-s} \binom{r-i}{s}\binom{i+t}{t} = \binom{r+t +1}{s+t +1} \text{,  }$$
(See \cite[3(b), p. 8]{Rio}.). \qed
\end{lemma}

The following result can be proved by straightforward induction.
\begin{lemma}\label{Summing}
For non-negative integers $n$ and $r$,
$$\sum_{j=r}^n \binom{j}{r}= \binom{n+1}{r+1} \text{.} $$ \qed 
\end{lemma}

\section{Combinatorial results for  \texorpdfstring{$\mathcal{OCI}_n $}{Lg}  }\label{CombOCIn}
Let  $A=\left\{(a_1,a_2,\ldots,a_m)| a_{i}\in \{1, 2, \ldots , n \} \text{ for } 1\leq i\leq m \right\}$ be a set of ordered tuples with fixed length $m$.
Define a relation $\mathbf{R}$ on the set $A$ as follows: For  $a=(a_1,a_2,\ldots,a_m)$ in $A$ and $b=(b_1,b_2,\ldots,b_m)$ 
in $A$,
$$a \ \mathbf{R}\ b \Leftrightarrow |a_{i}| \leq |b_{i}| \text{ for all } 1\leq i\leq m \text{.}$$

\begin{lemma} \label{cont}  A transformation $\alpha$ is a contraction if and only if  \\ $g(\im\alpha)\ \mathbf{R}\ g(\dom\alpha)$.  
\end{lemma} 
\begin{proof}
 The result follows from the definition of a contraction.
\end{proof}

For a tuple $\mathbf{r}=(r_1,r_2,\ldots, r_p)$ with $1\leq r_1<r_2<\cdots<r_p\leq n$, let 
$$ a_{\mathbf{r}}=\{\alpha\in \mathcal{OCI}_n | \im\alpha=\mathbf{r}\}  \text{.}$$
Let $g(\mathbf{r})=(d_1,d_2,\ldots,d_{p-1})$ and $\sum d_i=q$. Then
\begin{lemma}\label{SameIm} $$\left| a_{\mathbf{r}}\right|=\binom{n-q+p-1}{p}.$$\end{lemma}
\begin{proof} 

The number of maps in $\mathcal{OCI}_n$ with image $\mathbf{r}$ is equal to the number of  possible domain
sets which we obtain by  summing over \vspace{2 pt} the number of  possible expansions of the image set by inserting $i$ extra spaces, $\binom{i+p-2}{p-2}$, \vspace{2 pt} 
multiplied  by the number of possible  transversal shifts, $n-q-i$. 
\begin{eqnarray*}
\left| a_{\mathbf{r}}\right|&=& \sum_{i=0}^{n-1-q}(n-q-i)\binom{i+p-2}{p-2}\\
&=& \sum_{i=0}^{n-1-q}\binom{n-q-i}{1}\binom{i+p-2}{p-2}\\
&=& \binom{n-q+p-1}{p}    \hspace{64 pt} \text{(by Lemma~\ref{VarVander})}  \text{.}
\end{eqnarray*}  \end{proof}

\begin{corollary}\label{GapSum} Let $\mathbf{r}=(r_1,r_2,\ldots, r_p)$ with $1\leq r_1<r_2<\cdots<r_p\leq n$ and $g(\mathbf{r})=(d_1,d_2,\ldots,d_{p-1})$. 
Let  $\mathbf{s}=(s_1,s_2,\ldots, s_p)$ with $1\leq s_1<s_2<\cdots<s_p\leq n$ and $g(\mathbf{s})=(d_1',d_2',\ldots,d_{p-1}')$. 
If $\sum_{i=1}^{p-1} d_i=\sum_{i=1}^{p-1} d_i'$ then $\left| a_{\mathbf{r}}\right|=\left| a_{\mathbf{s}}\right|$.\end{corollary}
\begin{proof} The result follows from Lemma~\ref{SameIm}, since $\left| a_{\mathbf{r}}\right|$ depends only on $\sum d_i=q$ and $p$.  \end{proof}

\bigskip
       
\begin{theorem}\label{OneProduct}
Let $\displaystyle S =\mathcal{OCI}_n$. Then for $p\geq 1$,   $$F(n;p) = n \binom{n+p-1}{2p-1} +  (1-p) \binom{n+p}{2p} .$$
\end{theorem}
\begin{proof} From the results above, we can count $F(n;p)$ by summing over possible images and their possible domains as follows. 
\allowdisplaybreaks{
\begin{eqnarray*}
&& \hspace{-20 pt} F(n;p) =  \\[0.5 em] 
&=& \sum_{q=p-1}^{n-1}\binom{q-1}{p-2} (n-q)\binom{n-q+p-1}{p} \\[0.5 em]
&=&\sum_{i=0}^{n-p}\binom{p-2+i}{p-2} [(n - p + 1) - i]\binom{n-i}{p}  \hspace{62 pt} (i = q - p + 1)  \\[0.5 em]
&=& (n-p+1)\binom{n+p-1}{2p-1}-\sum_{i=0}^{n-p}i\binom{n-i}{p}\binom{p+i-2}{p-2}   \hspace{9pt} \text{(by Lemma~\ref{VarVander})} \\[0.5 em] 
&=& (n-p+1)\binom{n+p-1}{2p-1}-\sum_{i=1}^{n-p} (p-1)  \binom{n-i}{p} \frac{(p+i-2) \cdots (i+1)i}{(p-1)!}       \\[0.5 em] 
&=& (n-p+1)\binom{n+p-1}{2p-1}-\sum_{i=1}^{n-p}(p-1)\binom{n-i}{p}\binom{p+i-2}{p-1}\\[0.5 em]
&=&  (n-p+1)\binom{n+p-1}{2p-1}-(p-1)\binom{n+p-1}{2p}  \hspace{20 pt}  \text{(by Lemma~\ref{VarVander})}  \\[0.5 em] 
&=&  n \binom{n+p-1}{2p-1} + (1-p) \binom{n+p-1}{2p-1} + (1-p) \binom{n+p-1}{2p} \\[0.5 em]
&=&  n \binom{n+p-1}{2p-1} +  (1-p) \binom{n+p}{2p} \text{.}
\end{eqnarray*}}
\end{proof}

The height, $h(\alpha)$, the right and left shoulder, $\varpi^+(\alpha)$ and $\varpi^-(\alpha)$,
the right and left waist, $w^+(\alpha)$ and $w^-(\alpha)$, and the fix, $f(\alpha)$, are
defined in Section~\ref{introduction}. 
To compute $F(n;p,m)$, the  number of maps in $\mathcal{OCI}_n$ with height $p$ and $m$ fixed points, we introduce further notation.

Let $f^-(\alpha)= \min \{x \in \dom \alpha : x\alpha =x \}$ and  $f^+(\alpha)= \max \{x \in \dom \alpha : x\alpha =x \}$.

Let $h^-(\alpha)= | \dom \alpha \cap \{ 1,2, \ldots ,  f^-(\alpha)-1 \} |$ and  $h^+(\alpha)= |  \dom \alpha \cap \{ f^+(\alpha)+1,f^+(\alpha)+2, \ldots ,n \} |$.

\begin{lemma} Let $\alpha \in \mathcal{OCI}_n$.
The set of fixed points of $\alpha$  is convex with respect to $\dom \alpha $. 
That is, if $x\in \dom \alpha $ such that $f^-(\alpha) \leq x \leq f^+(\alpha)$, then 
$x \alpha = x$.   
\end{lemma}
\begin{proof}
 The result follows from the contraction property, as shown by \\  Adeshola and Umar  \cite[Lemma 1.1]{Ades}.
\end{proof}

For a semigroup $S$, we define 
\begin{align*}
& F(n;l^-,l^+,\lambda^-, \lambda^+, m,m^-,m^+,p,p^-,p^+)=  \\[0.5 em]
  &  =|\{ \alpha \in S :   (\varpi^-(\alpha)=l^- )\wedge  (\varpi^+(\alpha)=l^+ )
\wedge (w^-(\alpha)= \lambda^-)  \wedge \\[0.5 em]  &  \hspace{20 pt}  \wedge  (w^+(\alpha)= \lambda^+)
  \wedge  (f(\alpha)=m)  
 \wedge  (f^-(\alpha)=m^-) \wedge  
  (f^+(\alpha)=m^+) \wedge  \\[0.5 em]
  &   \hspace{20 pt} \wedge (h(\alpha)=p )  \wedge (h^-(\alpha)=p^- ) 	\wedge (h^+(\alpha)=p^+ )\}  |.
\end{align*}

\begin{theorem}\label{TwoProduct}  Let $\displaystyle S =\mathcal{OCI}_n$.  If $m=p$, then $F(n;p,m)=\binom{n}{m}$. If $m < p$, then 
 $$F(n;p,m)= (p-m-1)\binom{n+p-m-2}{2p-m } + 2 \binom{n+p-m-1}{2 p-m} \text{.}$$
\end{theorem}
\begin{proof}
The case $m=p$ is clear. For $m<p$ and $m\neq 0$,
as in the proof of Lemma~\ref{SameIm}, 
the maps can be counted by considering the number of possibilities of obtaining a domain set by expanding the image set.
Due to the contraction condition and since there are fixed points, right or left shifts of the domain set are not possible.

Along the lines of Al-Kharousi et al.~\cite{Kha1},
a transformation $\alpha \in \mathcal{OCI}_n$ can be split into three parts, an increasing part, followed by a fixed part,
and then a decreasing part. 
If a mapping has all three parts, then the set of fixed points is nonempty and occurs in the middle. In this case the number of spaces available
for expansion on the left of $m^-$ is $\lambda^- - l^-$.
One of these spaces necessarily needs to be inserted just before $m^-$.
Thus effectively we need 
to take into account all the possibilities of separating $\lambda^- -l^- -1$ objects with $p^- -1$ separators.
These $p^- -1$ separators correspond to the elements in the image which are between the smallest one and the smallest of the fixed points.
We need to multiply the number of possibilities of distributing these $p^- -1$ image points by the number of possibilities 
of expanding the image set, and then take the sum over all possible $\lambda^-$. 
Using Lemma~\ref{VarVander}, we obtain
$$ \sum_{\lambda^-=l^-+1}^{m^- -p^-} \binom{m^- -\lambda^- -1}{p^- -1}  \binom{\lambda^- - l^- -1 +p^- -1}{p^- - 1}
=  \binom{ m^- -l^- + p^- -2 }{2p^- -1} \text{.}$$
Since $p^+= p-m-p^-$, a similar calculation allows one to determine 
the number of possibilities  of expanding the image set      on the right of the fixed points as
$$\binom{l^+-m^++(p-m-p^-)-2 }{2( p-m-p^-)-1} \text{.}$$
For given $m^-$ and $m^+$, the number of possibilities to distribute the fixed points is 
$$ \binom{m^+-m^- -1}{m-2} \text{.}$$

The total number of transformations in $\mathcal{OCI}_n$ with the three parts (an increasing, a fixed, and a decreasing part) is obtained by multiplying these three expressions and 
summing over all possible $l^+$, 
$l^-$,$m^+$, $m^-$ and $p^-$. Lemmas \ref{Summing} and \ref{VarVander} will be applied repeatedly. 
\allowdisplaybreaks{
\begin{align*}
& \ \ \ \ \sum_{p^-=1}^{p-m-1} \ \sum_{m^-=p^-+2}^{n-(p-p^-)} \ \sum_{m^+= m^-+m-1}^{n-(p-m-p^-)-1} \  \sum_{l^-=1}^{m^--p^--1} \ \sum_{l^+ =m^+ + (p-m-p^-)+1}^n  
\\[0.5 em]
& \ \ \ \ \binom{m^--l^-+p^--2}{2p^--1}\binom{l^+-m^++(p-m-p^-)-2}{2( p-m-p^-)-1} \binom{m^+-m^- -1}{m-2}\\[0.5 em]
& =\sum_{p^-=1}^{p-m-1} \  \sum_{m^-=p^-+2}^{n-(p-p^-)} \ \sum_{m^+= m^-+m-1}^{n-(p-m-p^-)-1}  \binom{m^+-m^- -1}{m-2}  \\[0.5 em]
& \ \ \ \ \sum_{l^-=1}^{m^--p^--1} \binom{m^--l^-+p^--2}{2p^--1} \\[0.5 em]
& \ \ \ \ \sum_{l^+ =m^+ + (p-m-p^-)+1}^n 
\binom{l^+-m^++(p-m-p^-)-2}{2( p-m-p^-)-1} \\[0.5 em]
& = \sum_{p^-=1}^{p-m-1} \ \sum_{m^-=p^-+2}^{n-(p-p^-)} \  \sum_{m^+= m^-+m-1}^{n-(p-m-p^-)-1} \\[0.5 em]
& \ \ \ \  \binom{m^+-m^- -1}{m-2}  \binom{m^--1+p^--1}{2p^-} 
\binom{n-m^++(p-m-p^-) -1}{2( p-m-p^-)} \\[0.5 em]
& = \sum_{p^-=1}^{p-m-1} \ \sum_{m^-=p^-+2}^{n-(p-p^-)} \binom{m^-+p^--2}{2p^-} \\[0.5 em]
& \ \ \ \ \sum_{m^+= m^-+m-1}^{n-(p-m-p^-)-1}    \binom{m^+-m^- -1}{m-2}
\binom{n-m^++(p-m-p^-) -1}{2( p-m-p^-)} \\[0.5 em]
&= \sum_{p^-=1}^{p-m-1} \ \sum_{m^-=p^-+2}^{n-(p-p^-)}\binom{m^-+p^--2}{2p^-} 
\binom{n-m^-+(p-m-p^-)-1}{2( p-p^-)-m-1	} \\[0.5 em]
&= \sum_{p^-=1}^{p-m-1} \ 
\binom{n+p-m-2}{2p-m } = (p-m-1)\binom{n+p-m-2}{2p-m }.
\end{align*}}

Now we consider the cases where the fixed points are at the beginning or at the end. Here we do the case where they are at the beginning, the other case is similar. 
In the case under consideration,  $p^-=0$. Using Lemmas \ref{Summing} and \ref{VarVander}, we
obtain the total number of transformations in $\mathcal{OCI}_n$ with a fixed part followed by a decreasing part. 
\allowdisplaybreaks{
\begin{align*}
 & \ \ \ \ \sum_{m^+ =m}^{n-(p-m)-1} \ \sum_{l^+ =m^+ + (p-m)+1}^n \binom{m^+ -1}{m -1} \binom{l^+-m^++(p-m)-2}{2( p-m)-1} \\[0.5 em]
 &=  \sum_{m^+ =m}^{n-(p-m)-1} \binom{m^+ -1}{m -1} \sum_{l^+ =m^+ + (p-m)+1}^n  \binom{l^+-m^++(p-m)-2}{2( p-m)-1} \\[0.5 em]
 &=  \sum_{m^+ =m}^{n-(p-m)-1} \binom{m^+ -1}{m -1}  \binom{n-m^++(p-m)-1}{2( p-m)}  =  \binom{n+p-m-1}{2 p-m}. 
 \end{align*}}
 
Finally, we need to consider the case where there are no fixed points. Again we do this by expanding the image set in order to obtain the domain set. 
If all points in 
the image which come before the $(i+1)$th point are shifted to the left by inserting $j+1$ spaces,  we get the following number of possibilities
$$\binom{ j + (i-1)}{i-1} \text{.}$$
If all the elements in the image which come after the $i$th point are shifted to the right
by inserting $k+1$ spaces, 
this  gives the following number of possibilities
$$\binom{ k + (p-i-1)}{p-i-1} \text{.}$$
Given $\lambda^-$ and $\lambda^+$, the number of possible images is 
$$\binom{\lambda^+- \lambda^- -1}{p-2} \text{.}$$
Summing over all possible $\lambda^-$ and $\lambda^+$, using  Lemmas \ref{Summing} and \ref{VarVander}, we obtain 
\allowdisplaybreaks{
\begin{align*}
 &  \ \ \ \ \sum_{i=1}^{p-1} \  \sum_{\lambda^+ = p+1}^{n-1} \  \sum_{\lambda^-=2}^{\lambda^+-p+1} \ \sum_{j=0}^{\lambda^- -2}   \ \sum_{k=0}^{n-\lambda^+ -1} \\[0.5 em] 
& \ \ \ \ \binom{\lambda^+- \lambda^- -1}{p-2} \binom{ j + (i-1)}{i-1} \binom{ k + (p-i-1)}{p-i-1} \\[0.5 em]
 &  = \sum_{i=1}^{p-1} \  \sum_{\lambda^+ = p+1}^{n-1} \ \sum_{\lambda^-=2}^{\lambda^+-p+1}      
\binom{\lambda^+- \lambda^- -1 }{p-2} \\[0.5 em]
&\ \ \ \  \sum_{j=0}^{\lambda^- -2} \binom{ j + (i-1)}{i-1} \sum_{k=0}^{n-\lambda^+ -1} \binom{ k + (p-i-1)}{p-i-1} \\[0.5 em]
 &= \sum_{i=1}^{p-1} \  \sum_{\lambda^+ = p+1}^{n-1} \ \sum_{\lambda^-=2}^{\lambda^+-p+1} 
\binom{\lambda^+- \lambda^- -1}{p-2} \binom{ \lambda^- + i-2}{i} \binom{ n -\lambda^+ + p-i-1}{p-i} \\[0.5 em]
&= \sum_{i=1}^{p-1} \  \sum_{\lambda^+ = p+1}^{n-1} 
\binom{\lambda^++i-2}{i+p-1}  \binom{ n -\lambda^+ + p-i-1}{p-i} \\[0.5 em]
&= \sum_{i=1}^{p-1} 
\binom{n+p-2}{2p} = (p-1) \binom{n+p-2}{2p}\text{.}
\end{align*}}
Here also there are additional cases to consider, namely when all elements in the image are shifted to the left to obtain the domain, or 
all elements in the image are shifted to the right to obtain the domain. We will do the case  where all elements
in the image are shifted to the right and double the number. Using  Lemmas \ref{Summing} and \ref{VarVander}, we obtain
\allowdisplaybreaks{
\begin{align*}
 &  \ \ \ \  \sum_{\lambda^+ = p}^{n-1} \  \sum_{\lambda^-=1}^{\lambda^+-p+1} \ \sum_{k=0}^{n-\lambda^+ -1} 
\binom{\lambda^+- \lambda^- -1}{p-2}  \binom{ k + p-1}{p-1} \\[0.5 em]
&=\sum_{\lambda^+ = p}^{n-1} \ \sum_{\lambda^-=1}^{\lambda^+-p+1} 
\binom{\lambda^+- \lambda^- -1}{p-2}  \sum_{k=0}^{n-\lambda^+ -1}  \binom{ k + p-1}{p-1} \\[0.5 em]
 &=  \sum_{\lambda^+ = p}^{n-1} 
\binom{\lambda^+ -1}{p-1}  \binom{n-\lambda^+ +p-1}{p}
= \binom{n+p-1}{2p} \text{.}
\end{align*}}
When $p=1$ and $m=0$, then it is clear that $F(n;p,m)= n(n-1)$.

The result of the theorem is obtained by taking the sum of all cases.
\end{proof} 

\begin{remark}
$F(n;p)\, $ of $\, \mathcal{OCI}_n\,$  can be obtained from  the 
expression $\,F(n;p,m)\,$ in Theorem~\ref{TwoProduct} by summing over all possible values of $m$.
\end{remark}

It turns out that the order of $\mathcal{OCI}_n$ can be related to Fibonacci numbers. We need the following proposition to show the relation.
\begin{proposition}\label{EvenOddFib} The alternating Fibonacci numbers can be evaluated by the following formulae.
\begin{itemize}
 \item[(i)]   $$F_{2n+1} = \sum_{p\geq 0} \binom{n+p}{2 p}\text{.}$$ This is Sequence \seqnum{A001519} in The On-Line Encyclopedia of
Integer Sequences~\cite{oeis} and satisfies the recurrence relation 
 $$a_n= 3a_{n-1} - a_{n-2} \text{,}$$with $a_0=1=F_1$ and $a_1=2=F_3$. 
 \item[(ii)]  $$ F_{2n} =\displaystyle{ \sum_{p\geq 0}\binom{n+p-1}{2 p-1}} \text{.}$$ This is Sequence \seqnum{A001906} in The On-Line Encyclopedia of
Integer Sequences~\cite{oeis} and satisfies the recurrence relation 
 $$a_n= 3a_{n-1} - a_{n-2}\text{,}$$ with $a_0=0=F_0$ and $a_1=1=F_2$. 
\end{itemize}  \qed
\end{proposition}

The order of $\mathcal{OCI}_n$ can be determined using the formula for $F(n;p)$ obtained in Theorem~\ref{OneProduct}.
\begin{theorem}\label{Sequences}
The order of $\mathcal{OCI}_n$, as a function of $n$, is equal to $$h_n= \frac{3n-1}5 F_{2n}-\frac{n-5}5F_{2n+1}\text{.}$$
This is Sequence \seqnum{A094864} in The On-Line Encyclopedia of
Integer Sequences~\cite{oeis} and satisfies the recurrence relation $$h_n=6h_{n-1}-11h_{n-2}+6h_{n-3}-h_{n-4}\text{,}$$
with $h_0=1$, $h_1=2$, $h_2=6$, $h_3=18$. 
\end{theorem}
\begin{proof}
We use the formula in Theorem~\ref{OneProduct} and sum over all possible $p$. For $p=0$, the order is $1$, since the empty map is the only element
with height $0$.
\allowdisplaybreaks{
\begin{align*} 
h_n&=|\mathcal{OCI}_n| \\[0.5 em]   
&= 1+ \sum_{p=1}^n \left( n \binom{n+p-1}{2p-1} +  (1-p) \binom{n+p}{2p}   \right) \\[0.5 em] 
&= 1+ n \sum_{p=1}^n  \binom{n+p-1}{2p-1} +  \sum_{p=1}^n \binom{n+p}{2p}-\sum_{p=1}^n p \binom{n+p}{2p} \\[0.5 em]
&= 1+ n \sum_{p\geq 0}  \binom{n+p-1}{2p-1} + \left( \sum_{p\geq 0} \binom{n+p}{2p} -1\right)-\sum_{p=1}^n p \binom{n+p}{2p} \\[0.5 em]
&=  n \sum_{p\geq 0}  \binom{n+p-1}{2p-1} + \sum_{p\geq 0} \binom{n+p}{2p}-\sum_{p=1}^n p \binom{n+p}{2p} \\[0.5 em]
&=  n F_{2n}+ F_{2n+1} -\underbrace{ \sum_{p=1}^n p \binom{n+p}{2p}}_{b_n}\text{.}
\end{align*}}
We obtain a recurrence relation on the sequence $b_n$ as follows.
\allowdisplaybreaks{
\begin{align*} 
b_n&= \sum_{p=1}^n p \binom{n+p}{2p}
= \sum_{p=1}^n \frac{2p}2 \frac{(n+p)!}{(2p)!(n-p)!} \\[0.5 em] 
&= \sum_{p=1}^n  \frac{1}2 \frac{(n+p)!}{(2p-1)!(n-p)!}= \sum_{p=1}^n \frac{(n+p)}2 \frac{(n+p-1)!}{(2p-1)!(n-p)!} \\[0.5 em]
&= \frac{n}2 \sum_{p=1}^n \frac{(n+p-1)!}{(2p-1)!(n-p)!}  + \sum_{p=1}^n \frac{p}2 \frac{(n+p-1)!}{(2p-1)!(n-p)!} \\[0.5 em]
&= \frac{n}2 \sum_{p=1}^n \binom{n+p-1}{2p-1}  + \frac14 \sum_{p=1}^n (2p-1+1) \frac{(n+p-1)!}{(2p-1)!(n-p)!} \\[0.5 em]
&= \frac{n}2 F_{2n}  + \frac14 \sum_{p=1}^n (2p-1) \frac{(n+p-1)!}{(2p-1)!(n-p)!} 
+ \frac14 \sum_{p=1}^n  \frac{(n+p-1)!}{(2p-1)!(n-p)!}\\[0.5 em]
&= \frac{n}2 F_{2n}  + \frac14 \sum_{p=1}^n  \frac{(n+p-1)!}{(2p-2)!(n-p)!} 
+ \frac14 \sum_{p=1}^n  \binom{n+p-1}{2p-1} \\[0.5 em]
&= \frac{n}2 F_{2n}  + \frac14 \sum_{p=1}^n (n+p-1) \frac{(n+p-2)!}{(2p-2)!(n-p)!} 
+ \frac14   F_{2n} \\[0.5 em]
&= \frac{2n+1}4 F_{2n}  + \frac14 \sum_{p=1}^n (n+p-1)  \binom{n+p-2}{2p-2} \\[0.5 em]
&= \frac{2n+1}4 F_{2n}  + \frac{n}4 \sum_{p=1}^n   \binom{n+p-2}{2p-2} +
\frac14 \sum_{p=1}^n (p-1)  \binom{n+p-2}{2p-2} \\[0.5 em]
&= \frac{2n+1}4 F_{2n}  + \frac{n}4 \sum_{p=0}^{n-1}   \binom{n-1+p}{2p} +
\frac14 \sum_{p=0}^{n-1} p  \binom{n-1+p}{2p} \\[0.5 em]
&= \frac{2n+1}4 F_{2n}  + \frac{n}4 F_{2n-1} +
\frac14 b_{n-1}\text{.} 
\end{align*}}
This recurrence relation can be used to eliminate $b_n$ from the above formula for $h_n$ as follows.
\allowdisplaybreaks{
\begin{align*} 
 h_n&=n F_{2n}+ F_{2n+1} -b_n \\
 h_{n-1}&=(n-1) F_{2n-2}+ F_{2n-1} -b_{n-1} \text{.}
\end{align*}}
Thus we can express $h_n$ in terms of $h_{n-1}$,
and after some simplifications which merely involve using the 
property $F_n=F_{n-1}+F_{n-2}$, we obtain
$$h_n=\frac14 h_{n-1} + \frac{n+2}4 F_{2n}+ \frac12 F_{2n+1}  \text{.}$$ 
Since $h_0=1$, $h_1=2$, $h_2=6$, $h_3=18$, it is straightforward to verify that $h_n$ 
is Sequence \seqnum{A094864} in The On-Line Encyclopedia of
Integer Sequences~\cite{oeis} which was studied by 
Barcucci et al.~\cite{Barc} and Rinaldi and Rogers~\cite{Rin}. A closed formula
for this sequence is 
$$h_n= \frac{3n-1}5 F_{2n}-\frac{n-5}5F_{2n+1}\text{.}$$ 

\end{proof}

\section{Combinatorial results for \texorpdfstring{$\mathcal{ODCI}_n$}{Lg}}\label{CombODCIn} 

In this section, we will obtain results analogous to the ones obtained in Section~\ref{CombOCIn} for the semigroup $\mathcal{ODCI}_n$.
For a semigroup $S$, we define 
\allowdisplaybreaks{
\begin{align*}
 F(n;k^-,k^+,l^+,p) & = 
  |\{ \alpha \in S :  
  (w^-(\alpha)= k^-) \wedge  (w^+(\alpha)= k^+) \wedge  \\[0.5 em]
  & \hspace{123 pt} \wedge  (\varpi^+(\alpha)=l^+ )  \wedge (h(\alpha)=p ) \}|\text{.} 
\end{align*}}
\begin{lemma}\label{product} Let $S=\mathcal{ ODCI}_n$, then 
$$ F(n;k^-,k^+,l^+,p)=  \binom{l^+ - k^+ + p-1}{p-1   } \binom{k^+-k^--1}{p-2} .$$
\end{lemma}
\begin{proof}
The number of possible images is  $\binom{k^+-k^--1}{p-2}$. This number needs to be multiplied by the number of possible
pre-images, which depends on $k^+$ and $l^+$. Because of the decreasing  property, there is only 
one direction in which the image set can be expanded to obtain the domain set, 
namely to the right. \vspace{3 pt}  There are $l^+ - k^+ $ extra spaces and $p-1$ separators, 
so the number of possibilities to expand the  \vspace{3 pt}  domain set is $\displaystyle{\binom{l^+ - k^+ + p-1}{p-1   }}$.
\end{proof}
This allows us to find the number of maps in $\mathcal{ODCI}_n$ of height $p$.
\begin{theorem}\label{height} Let $S=\mathcal{ ODCI}_n $. Then 
$$ F(n;p)=  \binom{n+p}{2p} .$$
\end{theorem}
\begin{proof}
Using the expression from Lemma~\ref{product}, we obtain
\allowdisplaybreaks{
\begin{align*}
 F(n;k^+,l^+,p)&=  \sum_{k^-=1}^{k^+-p+1}\binom{l^+ - k^+ + p-1}{p-1   } \binom{k^+-k^--1}{p-2}  \\[0.5 em]  
 &= \binom{l^+ - k^+ + p-1}{p-1   } \sum_{k^-=1}^{k^+-p+1} \binom{k^+-k^--1}{p-2}   \\[0.5 em]  
   &= \binom{l^+ - k^+ + p-1}{p-1   }  \binom{k^+-1}{p-1} \hspace{ 20 pt} \text{(by Lemma~\ref{Summing})}  \text{.}      
\end{align*}}
Summing over all possible right waists and right shoulders, using Lemma~\ref{VarVander}, we obtain 
\allowdisplaybreaks{
\begin{align*}
F(n;p)&= \sum_{l^+=p}^n \ \sum_{k^+=p}^{l^+}  F(n;k^+,l^+,p)  \\[0.5 em] 
&= \sum_{l^+=p}^n \ \sum_{k^+=p}^{l^+}  \binom{l^+ - k^+ + p-1}{p-1   }  \binom{k^+-1}{p-1} \\[0.5 em]
&= \sum_{l^+=p}^n   \binom{l^+  + p-1}{2p-1   } = \binom{n+p}{2p}.
\end{align*}}
\end{proof}
We can extend this results to compute $F(n;p,m)$ as follows. 
\begin{lemma}\label{2ndProduct} Let $S=\mathcal{ ODCI}_n $. For $m < p$, 
$$ F(n;k^-,k^+,l^+,m,p)=  \binom{l^+-k^++p-m-2}{p-m-1} \binom{k^+-k^--1}{p-2} .$$
\end{lemma}
\begin{proof}
The  proof is analogous to  \vspace{3 pt}  the one of Lemma~\ref{product}.
If there are $m$ fixed points, because of the contraction and decreasing properties, these 
are the first $m$. 
There are $l^+ - k^+  $ extra spaces and $p-m-1$ separators.
One of these spaces is used to ensure that there are  no more than $m$ fixed points, i.e., it
is inserted \vspace{3 pt} after the last of the fixed points. 
Thus the number of possibilities to expand the domain  \vspace{3 pt}  set is $\displaystyle{\binom{l^+-k^++p-m-2}{p-m-1}}$.

\end{proof}
This allows us to find the number of maps in $\mathcal{ODCI}_n$ of height $p$ with $m$ fixed points.
\begin{theorem}\label{fixed} Let $S=\mathcal{ ODCI}_n $. If $m=p$, then $F(n;p,m)=\binom{n}{p}$. If $m < p$, then
$$ F(n;p,m)=  \binom{n+p-m-1}{2p-m} .$$
\end{theorem}
\begin{proof}
 The case $m=p$ is clear. For $m < p$, we repeat the procedure in the proof of Theorem~\ref{height}, using the expression 
 from Lemma~\ref{2ndProduct}. Note that because $m \neq p$, the case $k^+=l^+$ is excluded. Thus we get
 \allowdisplaybreaks{
\begin{align*}
 F(n;k^+,l^+,p,m)&=  \sum_{k^-=1}^{k^+-p+1}\binom{l^+-k^++p-m-2}{p-m-1}\binom{k^+-k^--1}{p-2}  \\[0.5 em]  
   &= \binom{l^+-k^++p-m-2}{p-m-1}  \binom{k^+-1}{p-1} \hspace{ 20 pt} \text{(by Lemma~\ref{Summing})}  \text{.}    
\end{align*}}
Then we sum over all possible right waists and right shoulders. 
Using  Lemmas \ref{VarVander} and \ref{Summing}, we obtain 
\allowdisplaybreaks{
\begin{align*}
F(n;p,m)&= \sum_{l^+=p+1}^n \ \sum_{k^+=p}^{l^+-1}  F(n;k^+,l^+,m,p)  \\[0.5 em] 
&= \sum_{l^+=p+1}^n \ \sum_{k^+=p}^{l^+-1}  \binom{l^+-k^++p-m-2}{p-m-1}   \binom{k^+-1}{p-1} \\[0.5 em]
&= \sum_{l^+=p+1}^n   \binom{l^++p-m-2}{2p-m-1} = \binom{n+p-m-1}{2p-m}.
\end{align*}}
\end{proof}

\begin{remark}
$F(n;p)\, $ of $\, \mathcal{ODCI}_n\,$ can  be obtained from  the 
expression $\,F(n;p,m)\,$ in Theorem~\ref{fixed} by  summing over all possible values of $m$.
\end{remark}

The order of $\mathcal{ODCI}_n$ can be determined using the formula for $F(n;p)$ obtained in Theorem~\ref{height}. It turns out that 
the order, as a function of $n$, can be expressed in terms of a single Fibonacci number. 
\begin{theorem}
$|\mathcal{ODCI}_n|=F_{2n+1}$, where $F_n$ is the $n$th Fibonacci number. 
\end{theorem}
\begin{proof}
The result follows from Theorem~\ref{height} and Proposition~\ref{EvenOddFib}.
\end{proof}

\section{Combinatorial results for \texorpdfstring{$\mathcal{ORCI}_n$}{Lg}}\label{CombORCIn} 
Let $\mathcal{OCI}_n^+$ be  the set of order-reversing contraction mappings of a finite chain  $X_n$, defined as
$$ \mathcal{OCI}^+_n = \{ \alpha \in \mathcal{CI}_n : \forall x,y \in 
\dom \alpha, \ x \leq y \implies x\alpha \geq y\alpha \} \text{.}$$

\begin{lemma} \label{bijec} There is a bijection between $\mathcal{ OCI}_n$ and $\mathcal{ OCI}^+_n$.
\end{lemma}
\begin{proof} Let $\alpha\in \mathcal{ OCI}_n$. If $h(\alpha)=1$, then $\alpha\in \mathcal{ OCI}^+_n$. 
Let $h(\alpha)>1$ and $$ \alpha= \left( \begin{array}{cccc}
           a_1&a_2&\cdots & a_p \\
          a_1\alpha& a_2\alpha &\cdots & a_p\alpha\\
          \end{array} \right) \text{,} $$ where $1\leq a_1 < a_2 < \cdots <a_p\leq n$. By the definition of  $\mathcal{ OCI}_n$,
          $$1\leq a_1\alpha < a_2\alpha < \cdots <a_p\alpha\leq n\text{.}$$ Let $g(\dom \alpha)=\left(t_1,t_2,\ldots,t_{p-1}\right) $,
          and $g(\im \alpha)=\left(d_1,d_2,\ldots,d_{p-1}\right)$. 
          Then we have  
          $$ \alpha= \left( \begin{array}{ccccc}
           a_1&a_1+t_1&\cdots & a_1+t_1+\cdots+t_{p-1} \\
          a_1\alpha& a_1\alpha+d_1 &\cdots & a_1\alpha+d_1+\cdots+d_{p-1}\\
          \end{array} \right) \text{,}$$ and by the contraction property,	 $\ d_i\leq t_i$ for all $i$. 
          
           We define a function $\theta: \mathcal{ OCI}_n\rightarrow\mathcal{ OCI}^+_n$, with $\theta(\alpha)=\alpha'$ given by 
       \begin{eqnarray*}
          \alpha'&=& \left( \begin{array}{cccc}  a_1 & a_1+t_{p-1} &\cdots & a_1+t_{p-1}+\cdots +t_1 \\
                                          a_p\alpha & a_p\alpha-d_{p-1} & \cdots & a_p\alpha -d_{p-1}-d_{p-2}-\cdots -d_1\\
          \end{array} \right)  \\[0.3 em] 
         & = &\left( \begin{array}{cccc}  a_1 & a_1+t_{p-1}  &\cdots  \\
                                          a_1\alpha +d_1+\cdots +d_{p-1} & a_1\alpha +d_1+\cdots +d_{p-2}&  \cdots \\ 
          \end{array} \right. \\[0.3 em] 
         &  & \hspace{182 pt}  \left. \begin{array}{cccc}    \cdots & a_1+t_{p-1}+\cdots +t_1 \\
                                            \cdots & a_1\alpha\\ 
          \end{array} \right) 
          \text{.}
          \end{eqnarray*} The function  $\theta$ is well defined, and since $d_i\geq 0$ for all $i$, $\alpha'$ is order reversing. 
          The gap of the domain of $\alpha'$ is equal to the reverse of 
          the gap of the domain of $\alpha$, 
          $$g(\dom \alpha')=\left(t_{p-1},t_{p-2},\cdots,t_1\right)=g(\dom \alpha)^R \text{.}$$        
          The gap of the image of $\alpha'$
          is equal to minus the reverse of the gap of the image of $\alpha$,  
          $$g(\im \alpha')=\left(-d_{p-1},-d_{p-2},\cdots,-d_1\right)
          =-g(\im \alpha)^R \text{.} $$
          Because for all $i$ with $1\leq i\leq p-1$, we have $1\leq d_i\leq t_i$, it follows that  $\alpha'\in \mathcal{ OCI}^+_n$.
          \end{proof}

  \begin{remark}\label{twice}
    Note that $\mathcal{ ORCI}_n=\mathcal{ OCI}_n\cup\mathcal{ OCI}^+_n$ and $\mathcal{ OCI}_n\cap\mathcal{ OCI}^+_n=\{\alpha\in \mathcal{ OCI}_n :\  h(\alpha)\leq 1\}$. 
  \end{remark}
  The following result follows from Theorem~\ref{OneProduct}, Lemma~\ref{bijec}, and Remark~\ref{twice}. 
 \begin{theorem}\label{LastProd}
   Let $\displaystyle S=\mathcal{ORCI}_n$. If $p=1$, then $F(n;p) = n^2$.
For $p>1$, 
$$F(n;p) = 2n \binom{n+p-1}{2p-1} +  (2-2p) \binom{n+p}{2p} .$$  \qed
\end{theorem}

We determine the number of maps with $m$ fixed points, for a given height $p$.
\begin{theorem}
  Let $\displaystyle S=\mathcal{ORCI}_n$. 
  If $m=p$, then $F(n;p,m)=\binom{n}{m}$. If $m=1 < p$, then
  $$F(n;p,m_1)= 2(p-2)\binom{n+p-3}{2p-1 } + 4 \binom{n+p-2}{2 p-1}- n  \text{.}$$
  If $1<m < p$, then
  $$F(n;p,m)= (p-m-1)\binom{n+p-m-2}{2p-m } + 2 \binom{n+p-m-1}{2 p-m} \text{.}$$
\end{theorem}
\begin{proof}
From Lemma~\ref{bijec} and Remark~\ref{twice}, if follows that the number of maps in $\mathcal{ OCI}^+_n$ 
with one fixed point is equal to  the number of maps with one fixed point in $\mathcal{ OCI}_n$. In order to  find the number
of maps with one fixed point in $\mathcal{ ORCI}_n$, we add the number for $\mathcal{ OCI}_n$ and the number for  
$\mathcal{ OCI}^+_n$, and subtract the number of partial identities of height one, as otherwise they would be counted twice. 
Using $m=1$ in Theorem~\ref{TwoProduct}, this gives 
$$F(n;p,1)= 2(p-2)\binom{n+p-3}{2p-1 } + 4 \binom{n+p-2}{2 p-1}- n  \text{.}$$
In $\mathcal{ OCI}^+_n$, there  are clearly  no maps with more than one fixed point. Thus for $m\geq 2$, the 
number of maps with $m$ fixed points in $\mathcal{ ORCI}_n$ is equal to the number of maps with $m$ fixed points in $\mathcal{ OCI}_n$.
From Theorem~\ref{TwoProduct}, we get 
$$F(n;p,m)= (p-m-1)\binom{n+p-m-2}{2p-m } + 2 \binom{n+p-m-1}{2 p-m} \text{.}$$
\end{proof}

We will calculate the order of  $\mathcal{ORCI}_n$ and find an expression for the order 
in terms of two consecutive Fibonacci numbers. In light of Remark~\ref{bijec}, we can use 
results from Section~\ref{CombOCIn} to significantly reduce the amount of calculations. 
\begin{theorem}\label{secondLast}
The order of $\mathcal{ORCI}_n$, as a function of $n$, is given by 
$$|\mathcal{ ORCI}_n|=   \frac{6n-2}5 F_{2n}-\frac{2n-10}5F_{2n+1}-1-n^2\text{.}$$
\end{theorem}
\begin{proof}
From Theorem~\ref{Sequences}, we know that the order \vspace{3 pt}  of $\mathcal{ OCI}_n$ is \\ $\displaystyle{\frac{3n-4}5 F_{2n}-\frac{4n-10}5F_{2n-1}}$.
\hspace{3 pt} This was obtained by summing over $F(n;p)$ for $p\geq 1$ and adding $F(n;p_0)=1$ for the empty map. In  $\mathcal{ ORCI}_n$,
according to Remark~\ref{bijec}, the case $p=1$ needs to be handled separately as well. Using the expression from 
Theorem~\ref{LastProd}, this gives 
\allowdisplaybreaks{
\begin{align*}
|\mathcal{ ORCI}_n|&= 1 + n^2+ \sum_{p=2}^{n}\left( 2n \binom{n+p-1}{2p-1} +  (2-2p) \binom{n+p}{2p}   \right) \\[0.5 em]
&= 1 + n^2+ \sum_{p=1}^{n}\left( 2n \binom{n+p-1}{2p-1} +  (2-2p) \binom{n+p}{2p}   \right)- 2n\binom{n}{1} \\[0.5 em]
&= 1 + n^2+ 2\left(|\mathcal{ OCI}_n| -1\right)- 2n^2 \\[0.5 em]
&= 1 + n^2+ 2\frac{3n-1}5 F_{2n}-2\frac{n-5}5F_{2n+1}-2- 2n^2 \\[0.5 em] 
&= \frac{6n-2}5 F_{2n}-\frac{2n-10}5F_{2n+1}-1-n^2 \text{.}
\end{align*}}
\end{proof}

\section{Acknowledgment}
The authors would like to thank the referee 
for helpful comments and suggestions. Financial support from The Reasearch Council of Oman grant RC/SCI/DOMS/13/01 is gratefully acknowledged

\nocite{*}
\bibliographystyle{amsplain}

\end{document}

%% file: cdmstdcmds.tex
\usepackage{graphicx}
\usepackage{epstopdf}
\usepackage{amsthm,amsmath,amssymb,amsxtra,amscd}
\usepackage{verbatim}
\usepackage{indent}
\usepackage{rotating}
\usepackage{multirow}
\usepackage{multicol}
\usepackage{url}
\usepackage{algorithmic}
\usepackage{xypic}

\newtheorem{theorem}{Theorem}[section]

\newtheorem{lemma}[theorem]{Lemma}
\newtheorem{corollary}[theorem]{Corollary}

\newtheorem{proposition}[theorem]{Proposition}

\newtheorem*{example*}{Example}

\newtheoremstyle{myexample}{3pt}{3pt}{\rmfamily}{}{\itshape}{:}{ }{\thmname{#1}\thmnumber{ #2}\thmnote{ (#3)}}
\theoremstyle{myexample}

\newtheoremstyle{myremark}{3pt}{3pt}{\rmfamily}{}{\itshape}{:}{ }{\thmname{#1}}
\theoremstyle{myremark}
\newtheorem{remark}[theorem]{Remark}
\newtheorem*{observation*}{Observation}

\newtheoremstyle{conjecture}{3pt}{3pt}{\itshape}{}{\bfseries}{.}{ }{\thmname{#1}\thmnote{ (#3)}}
\theoremstyle{conjecture}

\newtheorem*{question*}{Question}

\newtheorem{theorem*}{Theorem}

\numberwithin{equation}{section}

\newcounter{algorithm}
\renewcommand{\thealgorithm}{\thesection.\arabic{algorithm}}
